\documentclass[12pt]{article}
\usepackage{verbatim,amsmath,amssymb}
\usepackage[isolatin]{inputenc}
\usepackage{epsfig}
\input math1.def
\newtheorem{theor}{Theorem}[section]

\newtheorem{lemma}[theor]{Lemma}

\newtheorem{prop}[theor]{Proposition}

\def\proof{\goodbreak\noindent{\sc Proof. }\nobreak}

\def\endproof{\par\nobreak\hbox to \hsize{\hfil\vrule width 5pt height
5pt}\goodbreak\vskip 3pt}

\begin{document}

\title{A thin-thick Decomposition for Hardy Martingales}
\author{Paul F. X. M\"uller \thanks{Supported
by the Austrian Science foundation (FWF) Pr.Nr. P15907-N08.}  }
\date{August $21^{th}$ 2010}
\maketitle

\begin{abstract} 
We prove  thin-thick decompositions, 
for  the class of Hardy martingales
and thereby  strenghten its   square function
characterization. 
We apply the underlying method  to 
several classical martinale inequalities, 
for which we give new proofs . 

\end{abstract}

\tableofcontents

\section{Introduction}
 Let $\bT^\bN = \{(x_i)_{i=1}^\infty  \}$ denote 
the countable product of the torus $\bT=  \{ e^{i\theta} : \theta \in
[0, 2\pi [ \}, $
equipped with its normalized Haar measure $\bP. $
A natural filtration of $\s -$ algebras on  $\bT^\bN$
is given by the coordinate projections
$$ P_k :  \bT^\bN  \to  \bT^k, \quad (x_i)_{i=1}^\infty \to
(x_i)_{i=1}^k .$$
Define   $\cF_k$ to be  the  $\s -$ algebra on  $\bT^\bN$ generated 
by $P_k.$

Let  $F = (F_k) $ be  an $L^1(\bT^\bN)-$bounded  martingale on the filtered 
probability space $(\bT^\bN, (\cF_k) , \bP). $
Conditioned on  $\cF_{k-1}$  the martingale differnce $\Delta F_k =
F_k- F_{k-1}$  defines an element in the Lebesque space of integrable,
 function of vanishing mean $L_0^1(\bT) . $ 
By definition the martingale $F = (F_k) $ belongs to the class of
Hardy martingles,  if, 
conditioned on  $\cF_{k-1},$
$$
\Delta F_k
 =F_k- F_{k-1}
\quad\text{defines an element in the Hardy space $H_0^1(\bT)
  . $ }$$

 Hardy martingales, introduced by Garling \cite{gar2},  arise throughout 
 Complex and Functional Analysis. For instance  
in renorming problems for Banach spaces \cite{dgt, xu1},
vector valued Littlewood Paley Theory \cite{xu2}, 
embedding problems \cite{b1}, isomorphic classification 
problems \cite{b2, m1}, factorization problems \cite{p1},
similarity problems \cite{pisim}, boundary convergence of vector valued analytic functions
\cite{e1,g1,gmau, glm},
Jensen measures \cite{busch,asmo}.

As pointed out by Garling \cite{gar2}, two robustness properties of  Hardy martingales are
particularily important for their use in Analysis.
\begin{enumerate}
\item The class of  Hardy martingales is closed under martingale
  transforms.
\item  For Hardy martingales, their $L^1$ norm is 
determined by  square functions. There exist 
$c, C >0$ so that for any Hardy martingale  $F = (F_k) $  ,
\begin{equation}\label{k3} c \bE |F_n| \le \bE ( \sum_{k=1}^n |\Delta F_k|^2 )^{1/2}
\le C \bE |F_n| .\end{equation} 
\end{enumerate}

In the present paper we  strengthen
 the square function characterization \eqref{k3}
for   Hardy martingales. 
We prove that every  Hardy martingale $F = (F_k)_{k = 1 }^n $ 
can be  written as 
$$ F =  G +  B  $$
where 
  $G = (G_k)_{k = 1 }^n $
and  $B = (B_k)_{k = 1 }^n $ are again Hardy martingales 
so that 
\begin{equation}\label{i1}
 \bE ( \sum_{k = 1 }^n \bE_{k-1} |\Delta G_k |^2 )^{1/2}
+ \bE ( \sum_{k = 1 }^n  |\Delta B_k |)
\le C\bE | F_n | .
\end{equation}
and 
\begin{equation}\label{i2} |\Delta G_k | \le A_0  | F_{k-1} | , \quad\quad k\le n . 
\end{equation}
The estimate  \eqref{i1} implies of course   the right hand side of
the square function estimate 
\eqref{k3} since  the triangle inequality and the 
Burkholder-Gundy  martingale
inequality \cite{sia} give 
\[
\begin{aligned}
\bE (\sum_{k=1}^n| \Delta F_k |^2)^{1/2} &\le 
\bE ( \sum_{k = 1 }^n |\Delta G_k |^2 )^{1/2}
+ \bE ( \sum_{k = 1 }^n  |\Delta B_k |)\\
&\le 2\bE ( \sum_{k = 1 }^n \bE_{k-1} |\Delta G_k |^2 )^{1/2}
+ \bE ( \sum_{k = 1 }^n  |\Delta B_k |).
\end{aligned}
\]

The uniform previsible estimate \eqref{i2} should be compared with 
uniform previsible estimates appearing in the classical Davis and Garsia
inequality \cite[Chapters III and IV]{sia}. 
(See also Section~\ref{davisgarsia}.)
For general martingales the Davis decomposition \cite[Chapter III]{sia}  
guarantees only uniform estimates by  previsible {\it and 
 increasing} functionals such as $\max_{m\le k-1} |F_{m}|. $
Hence a routine application of the  Davis decomposition could yield only 
$$ |\Delta G_k | \le A_0 \max_{m\le k-1} |F_{m}| . $$ 

The present paper  exploits a basic and general iteration principle 
extracted from the work of J. Bourgain \cite{b1}. In its simplest form
it  yields a comparison theorem  between square functions as follows: 
Assume that $u_1, \dots, u_n$ and  $v_1, \dots, v_n$
are non- negative, integrable functions so that the following 
set of estimates hold true,
$$ \bE ( \sum_{m=1}^{k-1} u_m^2   + v_k^2)^{1/2} \le \bE (
\sum_{m=1}^{k} u_m^2  )^{1/2},
\quad\quad k \le n. $$
Then we have
$$ \bE ( \sum_{m=1}^{n} v_m^2   )^{1/2} \le 2 \bE (
\sum_{m=1}^{n} u_m^2  )^{1/2}.
$$

\paragraph{Acknowledgement:} 
It is my pleasure to thank S. Geiss, S. Kwapien and
M. Schmuckenschl\"ager
for very helpful and informative conversations during the preparation 
of this paper.
\section{Basic Iteration} 
In this section we review J. Bourgain's iteration method introduced in \cite{b1}.
It provides  upper estimates for the norm in 
$L^1(\ell^2) .$ By its generality the iteration method  can easily be adapted to
a variety of different situations. In this paper we apply it to obtain
proofs of  four different martingale inequalities. 

Consider first the  elementary  Lemma.
\begin{lemma} \label{numbers}
Let $ 0 \le s \le 1 , $
and $A, B \ge 0 . $ Then
\begin{equation}\label{n1}
 Bs \le s^2 A + ( A^2 + B^2 ) ^{1/2} - A  . 
\end{equation}
\end{lemma}
\proof
Since  $ 0 \le s \le 1 , $
 we have  $1 - s^2 < ( 1- s^2 ) ^{1/2} .$
Multiply by $A>0$, add $ Bs$ and use the Cauchy-Schwartz inequality.  This gives  
$$ 
A ( 1 -s^2 ) + B s 
\le A ( 1-s^2) ^{1/2} + B s 
\le ( A^2 + B ^2 )^{1/2} . 
$$
Subtracting $A ( 1 -s^2 ) $ gives \eqref{n1}.
\endproof

Let $(\O, \bP ) $ be  a probability space and write  $\bE $ to denote expectation in $(\O, \bP ). $

\begin{theor} \label{eloneit}
Let $n \in \bN . $ 
 Let   $  u_1,  \dots, u_n  \in  L^1(\O) ,$
and form the partial sums
 $$Z_k = \sum_{m= 1 }^k u_m , \quad\quad\text{ $1 \le  k \le n .$} $$
Assume that  $   v_1,  \dots, v_n,$ and  $  w_1 , \dots,  w_n$
be non-negative  in $ L^1(\O) ,$ so that 
the following estimates hold 
\begin{equation}\label{elonehypo}
\bE ( |Z_{k-1}|^2  +  v_k^2 )^{1/2}  +  \bE w_k 
\le \bE |Z_k|
\quad\text{ for  $1 \le  k \le n .$} \end{equation}
Then 
\begin{equation}\label{elsum}
\bE (\sum_{k = 1 }^n v_k ^2 )^{1/2} 
+ 
\bE \sum_{k=1}^n w_k
\le 2
(\bE |Z_n|)^{1/2}(\bE\max_{k\le n} |Z_k|)^{1/2} .
\end{equation}
\end{theor}
\proof 
Let $0\le \e \le 1 $ be defined by 
\begin{equation}\label{eleps}
 \e^2  = (\bE |Z_n|)(\bE\max_{k\le n} |Z_k|)^{-1} . \end{equation} 
 Choose next non negative $ s_k \in L^{\infty}$ so that 
$\sum_{k=1}^n s_k^2 \le \e^2 . $ 
Apply  Lemma~\ref{numbers} with 
$$A = |Z_{k-1}|, \,  B =  v_k \,\,\text{ and }\,\, s =  s_k . $$
This yields the point-wise estimates
 \begin{equation}\label{elpoint}
 v_k s_k  \le s_k^2 |Z_{k-1}|   +  (|Z_{k-1} |^2 + 
v_k^2)^{1/2}  -  |Z_{k-1}|  
\end{equation}
Integrating the  point-wise estimates \eqref{elpoint}  gives
 \begin{equation}\label{elintegrate}
\bE ( v_k s_k )
 \le \bE(  s_k^2|Z_{k-1}|)
+\bE (|Z_{k-1} |^2 +v_k^2 )^{1/2}
-   \bE|Z_{k-1} |^2 .
\end{equation}
Next apply the hypothesis \eqref{elonehypo}   to the central term 
$\bE (|Z_{k-1} |^2 +v_k^2 )^{1/2} $ appearing in the integrated
estimates \eqref{elintegrate}. This gives
 \begin{equation}\label{elgives}
\bE ( v_k s_k)
 \le \bE ( s_k^2 |Z_{k-1} |)
+ \bE |Z_{k} |
-   \bE |Z_{k-1} | -\bE w_k .
\end{equation}
Taking the sum over $ k \le n $ and exploiting   the telescoping nature
of the right hand side of \eqref{elgives}  yields,
\begin{equation}\label{eltele}
\begin{aligned} 
 \bE(  \sum_{k=1}^n v_k s_k )
+ \sum_{k=1}^n \bE w_k
&\le \bE |Z_n|
+\bE (\sum_{k=1}^n  s_k^2 |Z_{k-1}| )
 \\
&\le 
\bE |Z_n|
+\e^2 \bE \max_{k\le n } |Z_{k-1}| 
\end{aligned}
\end{equation}
Since \eqref{eltele}   holds for  every choice of 
 $ s_k \in L^{\infty}$ such  that 
$ \sum_{k=1}^n s_k^2  \le \e^2 , $ we may take  the supremum 
and obtain, by duality,  
the square function estimate
$$ \e \bE ( \sum_{k =1 }^n v_k^2 )^{1/2} 
+  \sum_{k=1}^n \bE w_k
\le 
\bE |Z_n|
+\e^2 \bE\max_{k\le n } |Z_{k-1}|
. $$
It remains to divide by $0\le\e \le 1$ and take into account \eqref{eleps}. This gives
\begin{equation}\label{elend}
\begin{aligned} 
 \bE(  \sum_{k=1}^n v_k^2 )^{1/2}
+ \sum_{k=1}^n \bE w_k
&\le 
\e^{-1}\bE |Z_n|
+\e \bE \max_{k\le n } |Z_{k-1}| \\
&= 2(\bE |Z_n|)^{1/2}(\bE\max_{k\le n} |Z_k|)^{1/2} .
\end{aligned}
\end{equation}

\endproof 
Theorem~\ref{eloneit} gives  estimates between plain integrals; in particular 
martingale structures are neither part of its hypothesis nor of its conclusion.
Nevertheless in Section~\ref{hardy} we employ  Theorem~\ref{eloneit} 
to prove 
 an inequality  for Hardy martingales. We use it to
estimate the $L^1$ norm of perturbed  square functions by the 
$L^1$ norm of the martingale itself.  
In Section~\ref{app} we discuss classical martingale inequalties
involving
different forms of  square
functions. There we will use a version of Theorem~\ref{eloneit}
 that is suitably adapted to estimating quadratic
expressions.

\section{Decomposing  Hardy Martingales}\label{hardy}
In this section we state and prove the main theorems of this paper. 
In the first paragraphs we collect  probabilistic results 
used later in the proof. We record here  a stochastic proof of Bourgain's complex
convexity inequality. This underlines the probabilistic nature of 
 Theorem \ref{elbrown}. 
\paragraph{Hardy Spaces, Brownian Motion and Complex Convexity.}

Let $\bT = \{ e^{i\theta} : \theta \in [0, 2\pi [ \} ,$
equipped with its  normalized Haar measure $dm .$  
For $ h \in L^p (\bT ) $ we say that $h$ belongs to
the Hardy space  $H^p (\bT ) $ if the harmonic extension of $h$ to the
unit disk is analytic. If moreover $\int_\bT h dm  = 0 $
we write $h \in  H^p_0 (\bT ).  $  
Recall J. Bourgain's complex convexity inequality \cite{b1}: There exists $\a_0 >
0$ so that   
\begin{equation} \label{jbi}
\int_{\bT}( |z|^2 + \a_0^2
|h|^2)^{1/2}dm \le \int_{\bT} |z +
h|dm,  \quad\quad z \in \bC,\quad h \in H^1_0(\bT) . 
\end{equation}
M. Schmuckenschl\"ager informed me that the Bourgain's  proof \cite{b1}
of \eqref{jbi}
gives     $\a_0^2 = 1/27  .$

Let $(B_t)$ denote complex 2D-Brownian motion on Wiener space, 
 and $((\cF_t), \bP ),$ the associated filtered 
probability space. Put
$$ \t = \inf\{ t> 0 : |B_t| > 1 \}. $$
With the  following proposition we verify  Bourgain's complex convexity
inequality. The proof uses Ito's formula and  the 
Cauchy-Riemann
equations.

\begin{prop}\label{pmi2} Let $z \in \bC,$ $h \in H^1_0(\bT) $
and let $\rho \le \t $ be a stopping time.
Then for $\a ^2 \le 1/6 $    
\begin{equation} \label{pmi}
\bE( |z|^2 + \a^2|h(B_\rho)|^2)^{1/2} \le \bE |z +
h(B_\rho)| .
\end{equation}
\end{prop}
\proof  We may put $|z| = 1 . $ By Ito's formula \cite{durr} we have
the identities
$$ \bE( 1 + \a^2|h(B_\rho)|^2)^{1/2}  = 1 +\frac12 
\bE\int_0^\rho \Delta(( 1 + \a^2|h(B_s)|^2)^{1/2})ds ,$$
and 
$$ \bE( |1  +  h(B_\rho)|)  = 1 +\frac12 
\bE\int_0^\rho \Delta( |1 + h(B_s)| )ds .$$
We calculate the Laplacians and evaluate the integrands 
on the right hand side as follows
$$ \Delta( |1 + h(z)| )= \frac{|h'(z)|^2}{ |1 + h(z)| } ,$$
and 
$$
\Delta(( 1 + \a^2|h(z)|^2)^{1/2})= \a^2 |h'(z)|^2
( 2 +
\a^2|h(z)|^2)( 1 +
\a^2|h(z)|^2)^{-3/2}
$$
An elementary calculation shows that for $\a^2 \le 1/6, $
$$ \a^2 |1 + w| (2 + \a^2|w|^2) \le  ( 1 + \a^2 |w|^2)^{3/2} , \quad\quad w \in \bC .$$
Hence for $\a^2 \le 1/6 , $
\begin{equation}\label{elca}
 \frac{\a^2  ( 2 +
\a^2|h(B_s)|^2)       }{( 1 + \a^2|h(B_s)|^2)^{3/2}}\le
\frac{1}{ |1 + h(B_s)| }  .
\end{equation}
 Multiplying both sides of \eqref{elca} by  $|h'(B_s)|^2$
and integrating gives
$$\bE\int_0^\rho \Delta(( 1 + \a^2|h(B_s)|^2)^{1/2})ds 
\le 
\bE\int_0^\rho \Delta( |1 + h(B_s)| )ds .$$
\endproof
\paragraph{Remarks.} 
\begin{enumerate}
\item The above proof applies, 
mutatis mutandis, to Conformal Martingales on Wiener Space.
 Let $ X ,  Y$ be real-valued and  integrable  
on  Wiener space, 
 $((\cF_t), \bP ).$ Assume $X, Y$ have identical 
quadratic variation,and vanishing  co-variance process,
 $$\la X \ra _t - \la Y \ra _t = \la X , Y \ra _t = 0 ,\quad\quad\, t \ge
0 .$$
The proof of Proposition~\ref{pmi2} shows that for 
$Z = X+ i Y ,$ $  w \in \bZ $ and $\a^2 < 1/6 ,$  
\begin{equation} \label{pmi22}
\bE( |w|^2 + \a^2|Z|^2)^{1/2} \le \bE |w +
Z| .
\end{equation}
\item  A short analytic  proof of 
\eqref{jbi} was obtained by M. Schmuckenschl\"ager 
who based his agrument  on Green's identity in the
following form:
$$ \pa_r \int_0^{2\pi} \vp(re^{it}) dt = r^{-1} \iint_{D_r} \Delta \vp
(z) dA(z) , $$ where 
$D_r = \{z \in \cC : |z| \le r \} ,$ 
$\Delta $ denotes the Laplacian and $dA(z) $ the area measure.  
Thus Green's formula replaces the use of Brownian Motion 
Ito's lemma.
\end{enumerate}
   
\paragraph{Hardy Martingales.}
 Let $\bT^\bN = \{(x_i)_{i=1}^\infty  \}$ denote 
the countable product of the torus $\bT $
equipped with its product Haar measure.

Let $n \in \bN , $ and denote by $\cF_n $
the the $\s-$algebra on $\bT^\bN$ generated by the 
cylinder sets 
$$ \{(A_1, \dots, A_n , \bT, \dots , \bT,\dots ,)\},$$
where $A_i,\, i \le n $ are measurable subsets 
of $\bT . $ 
Thus $(\cF_n)$  is an increasing sequence of 
$\s-$algebras canonically linked to the product structure
of $\bT^\bN .$ 
Subsequently we let $\bE_{n}$ denote the conditional expectation with
respect to the $\s-$algebra  $\cF_n .$
Let $F = (F_n)$ be an  $(\cF_n)$ martingale in  $L^1(\bT^\bN ). $
Denote its difference sequence by 
$$ \Delta F_n = F_n - F_{n-1} . $$
By definition  $F = (F_n)$ is a  Hardy martingale
if for almost all $(x_1, \dots, x_{n-1}) \in\bT^{n-1} $
fixed, the function
$$y \to  \Delta F_n(x_1, \dots, x_{n-1},  y),$$
defines an element in $H^1_0 (\bT ).  $ 

As shown by J. Bourgain \cite{b1}, the complex convexity 
inequality \eqref{jbi} combined with
Theorem~\ref{eloneit} yields the following 
 square function estimate  for Hardy martingales 
\begin{equation}\label{sfe}
\bE ( \sum_{k =1}^n |\Delta F_k|^2 )^{1/2} \le 2 \a_0^{-1}(\bE |F_n|)^{1/2}( \bE\max_{k\le n} |F_k|)^{1/2}.
\end{equation}
The  
 B. Davis martingale 
inequality \cite[p.37]{sia}
\begin{equation}\label{bdi}
\bE\max_{k\le n} |F_k| \le \sqrt{10} \bE ( \sum_{k =1}^n |\Delta F_k|^2
)^{1/2} ,
\end{equation}
and \eqref{sfe} imply that 
\begin{equation}\label{sfc}
\bE ( \sum_{k =1}^n |\Delta F_k|^2 )^{1/2} \le C_0 \bE |F_n| ,
\end{equation}
with $C_0 = 4 \times \a_0^{-2} \times \sqrt{10} . $
The converse  is a consequence of 
\eqref{bdi}. With a different constant and by a different method, 
the estimate \eqref{sfc} was abtained by 
B. Garling ~\cite{gar2}.

We next state the main theorem of the present paper. 
It provides a thin-thick decomposition for Hardy martingales 
and 
 strengthens the square function characterization
\eqref{sfc}.  
\begin{theor}\label{hardydavis}
Let   $A_0 = 4 \a_0^{-1}  $ and  $C_1 = 4  \a_0^{-1}\sqrt{ 10}\times A_0 . $
Every Hardy martingale $F = (F_k)_{k = 1 }^n $
can be decomposed as 
\begin{equation}\label{b1}
F = G + B 
\end{equation} 
where $G = (G_k)_{k = 1 }^n $
and  $B = (B_k)_{k = 1 }^n $ are again Hardy martingales so that the
following holds:
\begin{enumerate} 
\item Integral bounds:
\begin{equation}\label{b3}
 \bE ( \sum_{k = 1 }^n \bE_{k-1} |\Delta G_k |^2 )^{1/2}
+ \bE ( \sum_{k = 1 }^n  |\Delta B_k |)
\le C_1\bE | F_n | .
\end{equation}
\item Previsible uniform estimates:
\begin{equation}\label{b33} |\Delta G_k | \le A_0  | F_{k-1} | , \quad\quad k\le n . 
\end{equation}\end{enumerate}\end{theor}
\paragraph{Comments.}
We emphasize several points in which the decomposition for 
Hardy martingales given above 
is distinct from  the classical Davis and Garsia inequality
\cite[Theorems III.3.5 and IV.4.3] {sia}, holding for
general martingales. 
We refer also to  Section~\ref{app} where we give an alternative proof of the classical Davis and Garsia
inequality  based  on the iteration method.
\begin{enumerate}
\item The right hand side of \eqref{b3} involves just the $L^1$ norm
  $\bE | F_n | $ and not the square function  
$\bE ( \sum_{k =1}^n |\Delta F_k|^2 )^{1/2}. $ 
\item
The decomposition \eqref{b1} of $  F_k $ as 
$$
  F_k =  G_k +  B _k , 
$$
yields {\it analytic} martingale differences $\Delta G_k$ 
and
$\Delta B_k$ in the following sense.  For $(x_1, \dots , x_{k-1} ) \in \bT^{k-1}$ fixed  
the martingale difference 
$$y \to   \Delta G_k(x_1, \dots , x_{k-1} , y) ,$$
defines an element in $H^\infty_0 (\bT ). $ 
Hence  the decomposing martingales   $G = (G_k)_{k = 1 }^n $ and   $B =
  (B_k)_{k = 1 }^n $
are in fact Hardy martingales.
\item The right hand side of the 
previsible estimate \eqref{b33} 
involves just the value of the martingale  at time $ {k-1}, $
and  {\it not}  the entire history of the martingale {\it up to time} 
 $ {k-1}. $ This fact reflects  J.Bourgain's complex
 convexity inequality \eqref{jbi} and, apparently,  does not follow from \eqref{sfc}.
\end{enumerate}

We obtain Theorem~\ref{hardydavis} from the basic iteration theorem
using as input the following thin-thick decomposition
in $H^1_0(\bT) . $

\begin{theor}\label{elbrown} Let $A_0 = 4  \a_0^{-1} . $
To  $h \in H^1_0( \bT) $  and $ z \in \bC  ,$
put 
$$ \rho = \inf \{ t< \t : |h(B_t)| > 2 \a_0^{-1}|z| \} , 
\quad\text{and}\quad 
 g(e^{i\theta}) = \bE ( h(B_\rho) | B_\t = e^{i\theta}) . $$
Then  
$ g  \in H^\infty_0( \bT) $ satisfies the integral bounds
$$( |z|^2 + A_0^{-2} \int_{\bT} |g|^2 dm ) ^{1/2} + A_0^{-1}
\int_{\bT}|h-g|dm
\le \int_{\bT} |z +
h|dm 
, $$
and  the uniform estimate
$$ |g| \le A_0|z| . $$  
\end{theor}
\paragraph{Comments.}
\begin{enumerate}
\item
In the course of proving the decomposition Theorem~\ref{elbrown}
we use stopping time arguments on the holomorphic martingale 
$h(B_t) , t\le \t$ together with J. Bourgain's complex convexity
inequality \eqref{jbi}. 
\item We gave a proof of the complex convexity inequality \eqref{pmi}
using  Ito's formula and Cauchy Riemann equations. Hence  
the 
 thin-thick decomposition of  Theorem~\ref{elbrown}  has   natural counterparts  in  
 pure stochastic settings.
\end{enumerate}
In Lemma~\ref{brown} we  separate the stopping time argument from the rest of the
proof.
\nocite{pfxm1, jopfxm} 
We use below the following observation of Varopoulos.
Let $\rho \le \t$ be a stopping time and let $h \in H^1_0(\bT ) .$
Then taking the expectation of $h(B_\rho ) $ conditioned to 
$\{B_\t =e^{i\theta} \}$ gives again an element in   $  H^1_0(\bT ) .$
Thus 
$$\theta \to  \bE (h(B_\rho ) | B_\t =e^{i\theta} )$$
is in  $ H^1_0(\bT ) . $ See \cite{var1}.

\begin{lemma}\label{brown}
Let $h \in H^1_0( \bT) $  and $M > 0 .$
Put 
$$ \rho = \inf \{ t< \t : |h(B_t)| > 2M \} , 
\quad\text{and}\quad 
 g(e^{i\theta}) = \bE ( h(B_\rho) | B_\t = e^{i\theta}) . $$
Then  
$ g \in H^\infty_0( \bT) $ satisfies the integral estimate 
\begin{equation}\label{te2}( M^2 + \frac1{12} \int_{\bT} |g|^2 dm ) ^{1/2} + \frac14
\int_{\bT}|h-g|dm
\le \int_{\bT}( M^2 +
|h|^2)^{1/2}dm 
,\end{equation} 
and the uniform bound
 \begin{equation}\label{te3}|g| \le 2M .\end{equation}
\end{lemma}
\proof
By homogeneity assume $M = 1 . $ Put
$$ X =  h(B_\t) , \quad X_0 =  h(B_\rho),\quad X_1 = X - X_0 , $$
and 
$$ A = \{\rho < \t \}, \quad B = \O \sm A . $$
Note first that  $ \bE X_0 = \bE X_1 = 0 , $
and 
\begin{equation}\label{te4}
 |X_0| \le 2,\quad \quad \text{supp}( X_1) \sbe A , \quad \quad
\bE |X_1| \le 2 \bE (1_A |X|) . 
\end{equation}
 Moreover by inspection, 
$$ 
\bE (1_A |X_0|) \ge 2 \bP(A) . $$

We will prove next that 
$$\bE (1  + |X|^2)^{1/2} \ge (1 +\frac1{12}\bE |X_0|^2)^{1/2} + \frac14 \bE 
|X_1|  .
$$
To this end we consider separately the contribution of 
the sets $A$ and $B$ to $\bE (1  + |X|^2)^{1/2}$ and verify the following two estimates
$$\bE (1_A(1  + |X|^2)^{1/2}) \ge \bP (A)+\frac14 \bE 
 |X_1| 
$$
and 
$$
\bE (1_B(1  + |X|^2)^{1/2}) \ge \bP (B) + \frac1{12}\bE |X_0|^2.
$$
Let $\cF_{\rho}$ denote the stopping time $\s-$algebra generated by 
 $\rho . $ We may then rewrite  
$$X_0 = \bE(X| \cF_{\rho}). $$ 
For  $ \o \in A ,$ we have $|X_0(\o)| \ge 21_A (\o) . $ Recall that 
$A = \{\rho < \t \}, $ hence $A$ is  $\cF_{\rho}$ measurable, and  
\begin{equation}\label{te5}
 \bE ( 1_A |X| ) \ge\bE (| 1_A \bE(X| \cF_{\rho})| ) \ge   2 \bP(A) .\end{equation} 
Using \eqref{te5} and \eqref{te4} we get
\begin{equation}\label{teila}
\begin{aligned}
\bE (1_A(1  + |X|^2)^{1/2})   &\ge \bE ( 1_A |X| ) \\
&\ge \frac12 \bE (
1_A |X| ) + \bP (A)\ge  \frac14 \bE 
 |X_1|  + \bP (A).  
\end{aligned}
\end{equation}

Next for $\o\in B ,$ $ |X(\o)| \le 2 . $
Recall next the elementary estimate $(1 + x)^{1/2} \ge 1 +x/3$ for
$0\le x \le 1. $  Hence
$$ (1  + |X(\o)|^2)^{1/2} \ge (1  + \frac14
  |X(\o)|^2)^{1/2}
\ge
1 +\frac1{12}|X(\o)|^2, \quad\quad \o \in B . $$
Next take expectations  and use $\bE|X^2| \ge \bE|X_0^2|$ to obtain
\begin{equation}\label{teilb}
\bE (1_B(1  + |X|^2)^{1/2})  \ge \bE (1_B (  1
+\frac1{12}|X|^2))
\ge \bP (B) + \frac1{12}\bE |X_0|^2.
\end{equation}
Add the  estimates \eqref{teila}and \eqref{teilb}. 
This gives,
\begin{equation}\label{teilc}
\begin{aligned}
\bE (1  + |X|^2)^{1/2}
 &\ge \bP(B) +\frac1{12}\bE |X_0|^2 + \bP (A)+\frac14 \bE 
|X_1|   \\
&\ge (1 +\frac1{12}\bE |X_0|^2)^{1/2} + \frac14 \bE 
|X_1|  .
\end{aligned}
\end{equation}

Finally  we use the above bounds  for $X, X_1, X_0$ on Wiener space to get  estimates
for 
$$ g(e^{i\theta}) = \bE ( X_0 | B_\t = e^{i\theta}) .$$
By a well known  observation of Varopoulos  \cite{var1}
 $ g \in  H^1_0(\bT ) . $ 
Moreover, since $g $ is obtained by conditional expectation 
from $X_0$ we get 
$$|g| \le 2 , \quad 
\int_{\bT} |g|^2 
dm
\le \bE |X_0|^2 ,\quad\text{and}\quad
\int_{\bT} |h-g| 
dm
\le \bE |X_1| 
$$
Since $ X =  h(B_\t) ,$
$$ 
\int_{\bT}( 1 +
|h|^2)^{1/2}
dm  =\bE ( 1 + |X|^2)^{1/2}
. $$
Combining with \eqref{teilc} gives
\begin{equation}\label{teild}
\begin{aligned}
\int_{\bT} ( 1 +|h|^2)^{1/2} dm 
&\ge (1 +\frac1{12}\bE |X_0|^2)^{1/2} + \frac14 \bE 
|X_1| \\
&\ge 
(1 + \frac1{12} 
\int_{\bT} |g|^2 dm
) ^{1/2}
 + \frac14 
\int_{\bT}|h-g|dm
\end{aligned}
\end{equation}
\endproof
We next merge the conclusion  of Lemma~\ref{brown} with 
the complex convexity estimate \eqref{jbi}. 
\paragraph{Proof of Theorem~\ref{elbrown}.}
Apply Lemma~\ref{brown} to $ \a_0^{}h  $
and
$M = |z| . $
This gives $ g \in H^\infty{(\bT)}$ so that 
$$ |g| \le 2  \a_0^{-1}|z|, 
$$
and 
\begin{equation}\label{fast}
( |z|^2 +  12 \a_0^{2} \int_{\bT} |g|^2 dm ) ^{1/2} +
 4 \a_0^{}
\int_{\bT}|h-g|dm
\le \int_{\bT}( |z|^2 +
\a_0^{2}|h|^2)^{1/2}dm 
. 
\end{equation}
It remains to invoke  \eqref{jbi}, asserting that the right hand side
of \eqref{fast} is bounded by 
$$ \int_{\bT} |z +h|dm .
$$
\endproof
Finally we give the details of the proof  of Theorem~\ref{hardydavis}.
We   show how to apply 
Theorem~\ref{elbrown} to obtain the thin-thick decomposition for Hardy
martingales. 
\paragraph{Proof of Theorem~\ref{hardydavis}.}
Fix $k \le n  $ and  $(x_1, \dots , x_{k-1})\in \bT^{k-1} .$ By assumtion
the martingale difference 
$$ h(y) = \Delta F_k (x_1, \dots , x_{k-1}, y )$$
defines an element in $H^1_0(\bT) .$ 
Put $$z =  F_{k-1} (x_1, \dots , x_{k-1}),$$
 and 
apply Theorem~\ref{elbrown}
to $h$ and $z. $ This gives a decomposition
$$ h = g + b, $$
with $g \in H^{\infty} _0(\bT) , $ and  $b \in H^1_0(\bT) ,$
so that 
$$|g| \le A_0|z|, 
                 $$
and 
$$( |z|^2 + A_0^{-2} \int_{\bT}|g|^2 dm ) ^{1/2} 
+ A_0^{-1} \int_{\bT}|b |d m  \le \int_{\bT} |z +
h|dm, $$ 

Define next with $(x_1, \dots , x_{k-1} ) \in \bT^{k-1}$ 
fixed above
$$ \Delta G_k(x_1, \dots , x_{k-1} , y ) = g(y),
\quad\quad\text{and} \quad
\Delta B_k(x_1, \dots , x_{k-1} , y ) = b(y).$$
Then we get the identity $\Delta F_k =\Delta G_k +\Delta B_k  $
and the estimates 
$$|\Delta G_k| \le A_0|F_{k-1}|
$$
together with 
$$( |F_{k-1}|^2 + A_0^{-2}\bE_{k-1} |\Delta G_{k}|^2) ^{1/2} 
+ A_0^{-1} \bE_{k-1}|\Delta B_k|\le \bE_{k-1}  |F_{k}| . $$ 
Taking expectations 
gives 
$$\bE( |F_{k-1}|^2 + A_0^{-2}\bE_{k-1} |\Delta G_{k}|^2) ^{1/2} 
+ A_0^{-1} \bE|\Delta B_k|\le \bE  |F_{k}|
. $$
Now apply Theorem~\ref{eloneit}
with 
$$u_k =\Delta F_{k},\quad 
v_k = (\bE_{k-1} |\Delta G_{k}|^2) ^{1/2} /A_0 
\quad\quad\text{and}\quad\quad
w_k =  |\Delta B_k |/ A_0  , \quad k\le n .
$$ 
This gives the estimate
\begin{equation}\label{c2}
 \bE ( \sum_{k = 1 }^n \bE_{k-1} |\Delta G_k |^2 )^{1/2}
+ \bE ( \sum_{k = 1 }^n  |\Delta B_k |)
\le 2A_0(\bE   | F_n |)^{1/2} (\bE  \sup_{k \le n }  | F_k |)^{1/2}.
\end{equation}
Next use  
the inequalities  of B. Davis \eqref{bdi} and
Bourgain/Garling  \eqref{sfc}
\begin{equation}\label{c3}
\frac1{\sqrt{10}}\bE  \sup_{k \le n }  | F_k |\le 
\bE ( \sum_{k = 1 }^n  |\Delta F_k |^2 )^{1/2}
\le
C\bE   | F_n | ,\end{equation}
where $C  =4 \a_0^{-2} \times\sqrt{10} .$ Inserting \eqref{c3} into
\eqref{c2} gives 
\eqref{b33}.

\endproof

\section{Further Applications}\label{app}
We continue with applications  of the iteration principle
to classical martingale inequalities.
We  deduce the previsible projection theorem,
the comparison theorem between square functions and conditional square 
functions, and  prove the Davis and Garsia inequality.

We start with  a variant of Theorem~\ref{eloneit}
that is adapted to  bounding quadratic expressions.
Let  $ (\Omega, \bP) $ be probability space and  denote by 
$\bE $ the expectation in $( \O, \bP) .$ 
\begin{theor} \label{iteration}
Let $n \in \bN . $ 
 Let   $  u_1,  \dots, u_n $
be non-negative  in $ L^1(\O) ,$
and 
 $$M_k = (\sum_{m= 1 }^k u_m ^2  )^{1/2}\quad\text{ for  $1 \le  k \le n .$} $$
Assume that  $   v_1,  \dots, v_n,$ and  $  w_1 ,\dots,  w_n$
be non-negative  in $ L^1(\O) ,$ so that 
the following estimates hold 
\begin{equation}\label{hypo}
\bE ( M^2_{k-1}  +  v_k^2 )^{1/2}  +  \bE w_k 
\le \bE M_k
\quad\text{ for  $1 \le  k \le n .$} \end{equation}
Then 
\begin{equation}\label{sum}
\bE (\sum_{k = 1 }^n v_k ^2 )^{1/2} 
+ 
\bE \sum_{k=1}^n w_k
\le 2
\bE  (\sum_{k = 1 }^n u_k ^2 )^{1/2}  .
\end{equation}
\end{theor}
\proof 
 Choose non negative $ s_k \in L^{\infty}$ so that 
$\sum_{k=1}^n s_k^2 \le 1 . $ 
Apply  Lemma~\ref{numbers} with 
$$A = M_{k-1}, \,  B =  v_k \,\,\text{ and }\,\, s =  s_k . $$
This yields the pointwise estimates
 \begin{equation}\label{point}
 v_k s_k  \le s_k^2 M_{k-1}   +  (M_{k-1} ^2 + 
v_k^2)^{1/2}  -  M_{k-1}  
\end{equation}
Integrating the  point-wise estimates \eqref{point}  gives
 \begin{equation}\label{integrate}
\bE ( v_k s_k )
 \le \bE( s_k^2M_{k-1} )
+\bE (M_{k-1}^2 +v_k^2 )^{1/2}
-   \bE M_{k-1}.
\end{equation}
Next apply the hypothesis \eqref{hypo}   to the central term 
$\bE (M_{k-1}^2 +v_k^2 )^{1/2} $ appearing in the integrated
estimates \eqref{integrate}. This gives
 \begin{equation}\label{gives}
\bE ( v_k s_k)
 \le \bE (  s_k^2 M_{k-1}  )
+ \bE M_{k}
-   \bE M_{k-1} -\bE w_k .
\end{equation}
Taking the sum over $ k \le n $ and exploiting   the telescoping nature
of the right hand side of \eqref{gives}  yields,
\begin{equation}\label{tele}
\begin{aligned} 
 \bE(  \sum_{k=1}^n v_k s_k )
+ \sum_{k=1}^n \bE w_k
&\le \bE M_n
+\bE (\sum_{k=1}^n  s_k^2 M_{k-1}
)
 \\
&\le 
\bE M_n
+\bE \max_{k\le n } M_{k-1} 
\end{aligned}
\end{equation}
Since \eqref{tele}   holds for  every choice of 
 $ s_k \in L^{\infty}$ such  that 
$ \sum_{k=1}^n s_k^2  \le 1 , $ we may take  the supremum 
and obtain, by duality,  
the square function estimate
$$  \bE ( \sum_{k =1 }^n v_k^2 )^{1/2} 
+  \sum_{k=1}^n \bE w_k
\le 
\bE M_n
+\bE\max_{k\le n } M_{k-1}
. $$
It remains to use  that clearly $M_n =\max_{k\le n } M_{k}
. $
\endproof

\subsection{Previsible Projections}
Let  $ (\Omega, (\cF_k ) , \bP) ,$ be filterd probability space. 
Integration in  $(\O,\bP)$ is written as  $\bE .$
Conditional expectation with respect to  $\cF_k  $ is denoted by
$\bE_k .$
Let $ (F_k ) $ be a martingale 
 in  $ (\Omega, (\cF_k ) , \bP) ,$
 and $ \Delta F_k = F_k - F_{k-1} .$

We prove next square function estimates for the 
sequence of previsible projections $\bE_{k-1}( |\Delta F_k|). $
With  different methods the following result was obtained in 
\cite{ping, b1, descha} and \cite[Section 5.6]{kwa}. 
\begin{prop}\label{previsible}
$$
 \bE(  \sum_{k = 1 } ^n  \bE_{k-1}^2( |\Delta F_k|) )^{1/2}
\le 
2 \bE(  \sum_{k = 1 } ^n   |\Delta F_k|^2 )^{1/2} .
$$ 
\end{prop}
\proof
Let  $ u : [0,1] \to \bC $ be integrable and $M > 0 . $
We  verify  next
the following elementary inequality,
\begin{equation}\label{ttt}
 ( M^2 + (\int_0^1 |u(t)|dt )^2)^{1/2} \le \int_0^1 (M^2 + |u(t)|^2)^{1/2}dt 
. 
\end{equation}
By normalization we may choose $M=1 . $
Fix $ a, b \in \bR$ with $ a^2 + b ^2 = 1 $ so that 
$$ ( 1 + (\int_0^1 |u(t)|  dt )^2)^{1/2} = a + b\int_0^1 |u(t)| dt$$
Now estimate simply as follows 
$$ 
\begin{aligned}
a + b\int_0^1 |u(t)| dt  &= \int_0^1  a  + b |u(t)| dt\\
                    & \le      (a^2 + b ^2)^{1/2}\int_0^1(1 + |u(t)|^2)^{1/2}dt.
\end{aligned}
$$
Since  $ a^2 + b ^2 = 1 $ the estimate \eqref{ttt} is verified.

By the Theorem \ref{iteration} the proof of Proposition~\ref{previsible}
is now immediate.  
Fix $ k \le n ,$ and form the square function
$$
M_k = (\sum_{m= 1 }^k |\Delta F_m|^2  )^{1/2}.
$$
An immediate application of \eqref{ttt} is,
$$   ( M^2_{k-1}  + \bE_{k-1}( |\Delta F_k|^2) )^{1/2}   
\le \bE_{k-1}  M_k . $$
Taking expectations gives 
$$  \bE ( M^2_{k-1}  + \bE_{k-1}( |\Delta F_k|^2) )^{1/2}   
\le \bE  M_k . $$
Now apply Theorem~\ref{iteration} 
with 
$$ u_k = |\Delta F_k|,
\quad
v_k = (\bE_{k-1}( |\Delta F_k|^2) )^{1/2} 
\quad\text{and}\quad w_k =0,  $$
to get the conclusion.
\endproof
\subsection{Burkholder-Gundy Inequality}
We prove  the Burkholder-Gundy estimate, see \cite[Theorem
III.4.3]{sia}
or \cite[Section 5.6]{kwa},
comparing  the martingale  square function to
the conditioned square functions.

Let $ (F_k ) $ be an integrable martingale 
 in a filtered probability space  $ (\Omega, (\cF_k ) , \bP) ,$
 with  differences  $ \Delta F_k = F_k - F_{k-1} .$
\begin{prop}\label{condi}
$$
\bE(  \sum_{k = 1 } ^n   |\Delta F_k|^2 )^{1/2}
\le 2 \bE \left(  \sum_{k = 1 } ^n  \bE_{k-1}( |\Delta F_k|^2)
\right)^{1/2}
$$ 
\end{prop}
\proof
Let  $ u : \bT \to \bC $ 
integrable and fix $M >0 .$
 By Minkowski's inequality,
\begin{equation}\label{sss}
\int_0^1 ( M^2 +  |u(t)|^2 )^{1/2}dt \le ( M^2 +\int_0^1 |u(t)|^2 dt )^{1/2} . 
\end{equation}
By Theorem~\ref{iteration} we reduced Proposition~\ref{condi}
to \eqref{sss}. 
Indeed,   fix $k \le n $ and form the conditioned square function
$$
M_k = (\sum_{m= 1 }^k \bE_{m-1}(|\Delta F_m|^2)  )^{1/2}.
$$
By  \eqref{sss},
$$  \bE_{k-1} ( M^2_{k-1}  +  |\Delta F_k|^2 )^{1/2}   
\le \bE_{k-1}  M_k . $$
Taking expectations gives 
$$  \bE ( M^2_{k-1}  +  |\Delta F_k|^2 )^{1/2}   
\le \bE  M_k . $$
Use  Theorem~\ref{iteration} 
with 
$$ u_k =  (\bE_{k-1}(|\Delta F_m|^2)  )^{1/2}, \quad
v_k =  |\Delta F_k| ,\quad\text{and}\quad w_k =0.  $$

\endproof
\subsection{Davis and Garsia Inequality}\label{davisgarsia}
Let  $ (\Omega, (\cF_k ) , \bP) $  be a 
 filtered probability space.

The martingale transform techniques of  Garsia \cite[Theorem IV.1.2]{sia},
applied to the original Davis decomposition \cite[Theorem III.3.5]{sia} of a martingale
 $F = (F_k)_{k = 1 }^n $
 into a previsible part  $G = (G_k)_{k = 1 }^n $
and  $B = (B_k)_{k = 1 }^n $  gives the inequality 
\begin{equation}\label{converse}
 \bE ( \sum_{k = 1 }^n \bE_{k-1} |\Delta G_k |^2 )^{1/2}
+ \bE ( \sum_{k = 1 }^n  |\Delta B_k |)
\le C\bE (\sum_{k=1}^n| \Delta F_k |^2)^{1/2} .
\end{equation} 
Thus the inequality \eqref{converse} is a consequence of  
separate theorems due to Davis {\it and} Garsia respectively.

We proceed by giving  a new   proof of  \eqref{converse} 
based on Theorem \ref{iteration} and a martingale thin-thick decomposition. 
\begin{theor}\label{deco}
Every martingale $F = (F_k)_{k = 1 }^n $ in $L^1$ admits 
a decomposition 
as 
$$F = G + B $$
where  $G = (G_k)_{k = 1 }^n $
and  $B = (B_k)_{k = 1 }^n $ are martingales 
so that the following holds:
\begin{enumerate}
\item Integral bounds:
\begin{equation}\label{clb2}
 \bE ( \sum_{k = 1 }^n \bE_{k-1} |\Delta G_k |^2 )^{1/2}
+ \bE ( \sum_{k = 1 }^n  |\Delta B_k |)
\le 8\bE (\sum_{k=1}^n| \Delta F_k |^2)^{1/2} .
\end{equation}
\item Previsible uniform estimates:
\begin{equation}\label{clb3} |\Delta G_k |^2 \le 2 
\sum_{m = 1 } ^{k-1} |\Delta F_m |^2  , \quad\quad k\le n . 
\end{equation}\end{enumerate}\end{theor}
\paragraph{Comments.} Following are two remarks relating to the
inequality \eqref{clb2} and to the nature of the uniform previsible estimates\eqref{clb3}.
\begin{enumerate}
\item
The lower estimate \eqref{clb2} 
is sharp in the following sense.
{\it Any}  martingale decomposition of  $F = (F_k)_{k = 1 }^n $
as 
$$ F_k = G_k' + B_k' $$
gives rise to a reciprocal  upper estimate.
By the triangle inequality and the conditional square function estimate
Theorem~\ref{condi}, 
\[
\begin{aligned}
\bE (\sum_{k=1}^n| \Delta F_k |^2)^{1/2} &\le 
\bE ( \sum_{k = 1 }^n |\Delta G_k' |^2 )^{1/2}
+ \bE ( \sum_{k = 1 }^n  |\Delta B_k' |)\\
&\le 2\bE ( \sum_{k = 1 }^n \bE_{k-1} |\Delta G_k' |^2 )^{1/2}
+ \bE ( \sum_{k = 1 }^n  |\Delta B_k' |).
\end{aligned}
\]
\item
The right hand side of the previsible uniform 
estimates \eqref{clb3} 
depends not only on the value of the martingale of $F$ at time 
$k-1$ but also on its history up to time $k-1 .$ To wit 
\eqref{clb3} involves
$$ |\Delta F_1|, \dots, |\Delta F_{k-1}| . $$
This aspect  contrasts the uniform previsible 
estimates in the thin-thick decomposition for 
Hardy martingales \eqref{b33}.
\end{enumerate}
We prove  Theorem~\ref{deco} by feeding Theorem~\ref{iteration} with 
 a  
thin-thick decomposition for  $L^1(\O)$.
It uses just truncation  and is a   simplified version of
Lemma~\ref{brown}. 
\begin{lemma}\label{Kfun}
 To each $ h \in L^1(\O)$ satisfying $\bE h = 0 $
and $M > 0 $ put
\begin{equation}\label{trunc}
g = 1_D h - \bE (1_D h) \quad{where} \quad D = \{ |h | \le 2 M \} . 
\end{equation}
Then 
$$ |g| \le 2 M , 
$$
{and}
$$( M^2 + \frac1{12}  \bE |g|^2 ) ^{1/2} + \frac14 \bE|h-g| \le \bE ( M^2 +
|h|^2)^{1/2} .$$
\end{lemma}
\proof By rescaling we may put $M = 1 . $
Let $A$ denote  the complement of $D, $ thus  $A = \{ |h | > 2  \}.$
To $ x \in A ,$ we have $|h(x)| \ge 21_A (x) . $ Hence
$$ \bE ( 1_A |h| ) \ge 2 \bP(A) , $$
and 
\begin{equation}\label{parta}
\bE( 1_A(1  + |h|^2)^{1/2})   \ge \bE ( 1_A |h| ) \ge \frac12 \bE (
1_A |h| ) + \bP (A). 
\end{equation}

Next for $x\in D ,$ $ |h(x)| \le 2 . $
Hence
$$ (1  + |h(x)|^2)^{1/2} \ge (1  + \frac14 |h(x)|^2)^{1/2}\ge
1 +\frac1{12}|h(x)|^2, \quad\quad x \in D . $$
Taking expectations  gives
$$ \bE (1_D(1  + |h|^2)^{1/2})  \ge \bE (1_D (  1 +\frac1{12}|h|^2))
. $$
Next recall that in  \eqref{trunc} we defined $ g =  1_D h - \bE (1_D
h). $ 
Hence $ \bE (1_D|h|^2)
\ge  \bE |g|^2 , $ and 
\begin{equation}\label{partd}
\bE( 1_D(1  + |h|^2)^{1/2})  \ge \bP (D) + \frac1{12}\bE |g|^2.
\end{equation}

Adding  \eqref{parta} and \eqref{partd} gives
\begin{equation}\label{cc}
\begin{aligned}
\bE (1  + |h|^2)^{1/2} 
 &\ge \bP(D) +\frac1{12}\bE |g^2| + \bP (A)+\frac12 \bE (
1_A |h| )  .
\end{aligned}
\end{equation}
Note that $\bE h = 0 $ implies 
$\bE (1_A h) = -\bE (1_D h) $ and 
$$ h - g = h - 1_D h +\bE (1_D h)= 1_A h-\bE( 1_A h) . $$
Moreover 
\begin{equation}\label{aa}
2  \bE (
1_A |h| )  \ge \bE | 1_A h-\bE (1_A h) | . \end{equation}
Inserting \eqref{aa} into  \eqref{cc} gives the result
\begin{equation}\label{bb}
\begin{aligned}
\bE(1  + |h|^2)^{1/2}  &\ge \bP(A)+\bP(D) +\frac1{12}\bE |g^2| +\frac14 \bE 
|h- g|  \\
&\ge( 1 +\frac1{12}\bE |g|^2)^{1/2} +\frac14 \bE 
|h- g|. 
\end{aligned}
\end{equation} 
\endproof

\paragraph{Proof of Theorem~\ref{deco}.}
Let $k \le n  $ and put
 $$M_{k-1} = (\sum_{m =1}^{k-1}  |\Delta F_m |^2 )^{1/2} . $$
 Lemma~\ref{Kfun} gives  a decomposition of
$\Delta F_k $ as  follows. Put 
$$D_k = \{ |\Delta F_k | \le 2 M_{k-1} \}, 
\quad\text{and}\quad  \Delta G_k = 1_{D_k}\Delta F_k  - \bE_{k-1}( 1_{D_k}\Delta F_k) . $$
Define $\Delta B_k$
by the  decomposition
$$\Delta F_k =\Delta G_k +\Delta B_k .$$ 
Then  $ \Delta G_k , \Delta B_k $ are $\cF_k $ measurable, and  
$$  \bE_{k-1} (\Delta G_k ) = \bE_{k-1} (\Delta B_k ) = 0 . $$
By construction
$$| \Delta G_k | \le 4 M_{k-1}. 
$$
By  Lemma~\ref{Kfun} 
\begin{equation}\label{ff}
 ( M_{k-1}^2 + \frac1{12}\bE_{k-1} |\Delta G_k |^2 )^{1/2} 
+  \frac14 \bE_{k-1} ( |\Delta B_k |) \le \bE_{k-1}  M_k . 
\end{equation}
Take expectations of \eqref{ff} to obtain
\begin{equation}\label{gg}
\bE ( M_{k-1}^2 + \frac1{12}\bE_{k-1} |\Delta G_k |^2 )^{1/2} 
+  \frac14 \bE ( |\Delta B_k |) \le \bE M_k . 
\end{equation}
Now  apply Theorem~\ref{iteration} 
with
$$ u_k = |\Delta F_k|,\quad v_ k = (\bE_{k-1} |\Delta G_k
|^2 )^{1/2}/4 ,\quad\text{and}\quad   
w_k =  |\Delta B_k |/4 ,\quad\quad k\le n .$$ 
\endproof

\bibliographystyle{abbrv}
\bibliography{hardymartingales}

\paragraph{Author Address}:\\    
Department of Mathematics\\
J. Kepler Universit\"at Linz\\
A-4040 Linz\\
pfxm@bayou.uni-linz.ac.at

\end{document}